\documentclass[11pt]{amsart}
\usepackage{cite,thmtools}
\usepackage[capitalise]{cleveref}
\usepackage{tikz}
    \usetikzlibrary{positioning,arrows}
\usepackage{stix} 
\declaretheorem[style = plain]{theorem}   
\declaretheorem[style = plain,      sibling = theorem]{corollary,lemma,proposition}
\declaretheorem[style = definition, sibling = theorem]{definition, example}
\declaretheorem[style = remark]{remark,problem}

\def\C*{{\textsl{C*}-algebra}}
\def\Cs*{{sub-\textsl{C*}-algebra}}
\newcommand{\Bh}{{\ensuremath{B(\Hilbert)}} }
\newcommand{\bk}{{B\tensor  K}}

\newcommand{\constant}[1]{\delta_\infty(#1)}

\newcommand{\Cstar}[1][]{\textsl{C*}{#1}} 
\newcommand{\dist}{{\textup{dist}}}
\newcommand{\geqCu}{\ensuremath{\underset{\widetilde{\quad}}{\succ}}}
\newcommand{\her}{{\textup{her}}}
\newcommand{\Hilbert}{\ensuremath{\mathcal{H}}}
\newcommand{\hs}{{\mathcal  H}}

\newcommand{\Kasp}{{\bar\kappa}}
\newcommand{\Mbb}{{\Mb}/{B}}
\newcommand{\Mbkbk}{{\Mbk}/({B\tensor   K})}
\newcommand{\Mbk}{{\mathcal M} (B\tensor   K)}
\newcommand{\Mb}{{\mathcal M} (B)}
\newcommand{\Sinf}{\bar\delta_{\infty}}
\newcommand{\tensor}{\otimes}
\renewcommand{\O}[1]{{\ensuremath{{\mathcal O}_{#1}} }} 

\begin{document}




\title{Double relative commutants in coronas of separable \textsl{C*}-algebras}

\author{Dan Kucerovsky}
\author{Martin Mathieu}

\address{$^1$D.  Ku\v{c}erovsk\'{y} \\
              University of New Brunswick at Fredericton \newline%
Canada\qquad E3B 5A3} \email{dkucerov@unb.ca}
\address{$^2$ Martin Mathieu\\Mathematical Sciences Research Centre\\ Queen's University Belfast, Belfast\\
Northern Ireland BT7 1NN}


\begin{abstract}
We prove a double commutant theorem for separable subalgebras of a wide class of corona \textsl{C*}-algebras,  largely resolving
a problem posed by Pedersen~\cite{pedersen.corona}. Double commutant theorems originated with von Neumann,
whose seminal result evolved into an entire field now called von Neumann algebra theory. Voiculescu later proved a
\textsl{C*}-algebraic double commutant theorem for subalgebras of the Calkin algebra.
We prove a similar result for subalgebras of a much more general class of so-called corona \textsl{C*}-algebras.
\end{abstract}


\maketitle

\section{Introduction}

\noindent
Let $S$ be a subset of an algebra~$D$. Its \textit{relative commutant\/} $S'$ in $D$ is defined by
\[
S'=\{y\in D\mid xy=yx \text{ for all }x\in S\};
\]
that is, the centraliser of $S$ in~$D$. Clearly, $S'$ is always a subalgebra of $D$, being unital if $D$ is.
The double relative commutant, $S''$, is $S''=(S')'$. In the case when $D=B(\hs)$ for a  Hilbert space~$\hs$,
the adjective `relative' is customarily dropped.
In a similar vein, the \textit{unitisation\/} of a subalgebra $A$ of a unital algebra $D$ is the algebra generated
by $A$ and the identity of~$D$.

The most fundamental result in all of von Neumann algebra theory is arguably von Neumann's double commutant theorem,
published in 1929 (see~\cite{vonNeumann}). We phrase the theorem as follows:
\begin{theorem}
Given a *-subalgebra of $B({\hs})$, the double commutant of the subalgebra is equal to the weak* closure of
its unitisation.
\end{theorem}
Approximately half a century later, Voiculescu~\cite{voi0,voi} proved a \textsl{C*}-algebraic version of the above theorem:
\begin{theorem}
Let $C(\hs):=B({\hs})/K({\hs})$ be the Calkin algebra of a separable
infinite-dimensional Hilbert space~${\hs}$. The double relative commutant of a
separable \Cs* is the unitisation of that subalgebra.
\end{theorem}
Recall that the multiplier algebra $\Mb$ of a given \C* $B$
is the idealizer of $B$ in its enveloping von Neumann algebra~$B^{**}$.  Since the multiplier algebra of the
compact operators ${ K}:=K({\hs})$ is $B({\hs})$, we may reasonably regard
the corona algebras $\Mbb$ as a generalization of the Calkin algebra
considered by Voiculescu. At a conference in 1988, Pedersen posed the problem of generalizing Voiculescu's theorem
to the setting of general corona algebras~\cite{pedersen.corona}, and this note provides a partial answer
to his problem. We use some ideas from Kadison-Singer \cite{KadisonSinger}, a theorem from $KK$-theory \cite{EK}, and  previous double commutant theorems \cite{EKrelcomm,GN} in our proof.
In~\cite{EKrelcomm} is the following double relative commutant
theorem for hereditary subalgebras in \C*s of the form $\Mbb$ with $B$ simple and separable:

\begin{theorem}\label{th:main}
Let $B$ be a separable simple \C*.
Let $C$ be a hereditary $\sigma$-unital \Cs* of the corona algebra $\mathcal M(B)/B$.
Then the double relative commutant of $C$ in the given corona
algebra is equal to the unitisation of\/~$C$.
\end{theorem}

Giordano and Ng~\cite[Corollary~3.5]{GN} proved:

\begin{theorem}\label{th:gn}
Let $B$ be a stable separable \C* and suppose that either $B$ is the compact operators or $B$ is  simple and purely infinite.
Then separable unital \Cs*s of the corona are equal to their own double relative commutant.
\end{theorem}

Corollary 4.12 in \cite{hadwin} shows that Voiculescu's original double commutant theorem does not
generalize to the Calkin algebra of a nonseparable Hilbert space.

\section{Extensions of \textsl{C*}-algebras and absorption properties}
\noindent
Let $A$ and $B$ be \C*s, with $A$ unital, $B$ separable and simple.  An extension
\[
0\longrightarrow B \longrightarrow C \longrightarrow A \longrightarrow0
\]
will be said to be \textit{unital} if $C$ is unital. Recall that an extension of $B$ by $A$ is
determined up to strong isomorphism by its Busby map --- the naturally associated map from $A$
to the quotient multiplier algebra, or corona algebra, of $B$, $\Mbb$. If $B$ is stable,
so that the Cuntz
algebra $\O{2}$ may be embedded unitally in $\Mbb$, then there is a natural notion of addition of extensions.
Recall that an extension $0  \to  B  \to
C \to  A  \to  0$ of $B$ by $A$ is said to be \textit{trivial\/} if there exists a splitting homomorphism $\pi\!\!:\> A\to C.$
Let us further say that an extension of \Cstar-algebras $0  \to  B  \to
C \to  A  \to  0$ is \emph{trivial in the nuclear sense} if the
splitting homomorphism  may be chosen to be weakly
nuclear as defined by Kirchberg in \cite{KirchbergUnpublished}:  the splitting homomorphism
$\pi\!\!:\> A\to C$ will be said to be \emph{weakly nuclear} if, for every $b\in
B\subseteq C$, the map
$$
A \ni a \ \mapsto \ b\pi (a) b^* \in B\subseteq C
$$
is nuclear.

Kasparov \cite{KasparovSV,KasparovKK} introduced a property that he called {absorbing}. Letting $\tau_1$ and $\tau_2$ denote Busby maps of  extensions, and letting  $v_i$ denote isometries generating a copy of $\O{2},$ we will say that the extension $\tau_1$ is  \emph{absorbing in Kasparov's sense\/} if it is unitarily equivalent, by a multiplier unitary,
to its sum $a\mapsto v_1 \tau_1(a)v_1^* + v_2 \tau_2(a) v_2^*$
with any trivial extension $\tau_2.$ For technical reasons, Kasparov assumes at one key point that the algebra  $A$ is nuclear (see Theorem 1 of \cite[\S7]{KasparovKK}). Moreover, an extension that is absorbing in Kasparov's sense is not unital, and we will need to consider unital extensions. Thus we make the following adjustment to the definition, where following \cite{EK}, we define, with $\tau_i$ and $v_i$ as above:
\begin{definition} A unital extension $\tau_1$ is \emph{absorbing in the nuclear sense\/} if it is unitarily equivalent
to its sum $a\mapsto v_1 \tau_1(a)v_1^* + v_2 \tau_2(a) v_2^*$
with any unital  extension $\tau_2$ that is trivial in the nuclear sense.
\label{def:nuc.absorbing}\end{definition}
 { Busby maps take their values in the corona, so this sum is in the corona.} Different choices of unitally embedded copies of $\O{2}$ lead to different notions of sums of extensions that turn out to be equivalent under the equivalence relation given by unitary equivalence of extensions.  The main result of \cite{EK} is a \C*ic characterization of the property of being absorbing in the nuclear sense; this algebraic property is called the \textit{purely large\/} property:
\begin{definition} We say that an extension is \emph{purely large} if for every positive element $c$ of the extension algebra that is not contained in the canonical ideal $B,$ there exists a stable subalgebra $D \subset \overline {cBc}$ which is full in $B.$ A positive element $c$ is said to have the \emph{purely large property} if it is not in $B$ and if there exists a stable subalgebra $D \subset \overline {cBc}$ which is full in $B.$
\label{def:purely.large}\end{definition}

Our definition of  absorbing in the nuclear sense specifies that the extensions absorbed are unital. If one of the ambient \Cstar-algebras is nuclear, then the weak nuclearity condition is automatic, and in this case, the term unital absorption is often used. Then we may say that the next lemma characterizes unital absorption: see \cite[Section 17]{EK} and the correction for the nonunital case in \cite{gabe2016}. This particular statement of the lemma is from \cite[Remark 2.9]{GR2020}.

\begin{lemma}[\!\protect{\cite[Remark 2.9]{GR2020}}]\label{lem:inf.repeats}
Let $A$ and $B$ be separable \C*s, with
$B$ stable and nuclear.  Let $w_i$ denote the generators of a unital copy of $\O{2}$ in the multipliers.
Consider  a unital essential  \C* extension $\tau$ of $B$ by~$A$.
The following are equivalent:
\begin{enumerate}
\item[(i)] the extension $\tau$ is unitally absorbing, meaning that it is unitarily equivalent
to its sum $w_1 \tau (a)w_1^* + w_2 \sigma(a) w_2^*$
with any trivial unital  extension $\sigma,$  and
\item[(ii)] the extension algebra, in $\Mb$, is purely large.
\end{enumerate}
\end{lemma}
Following \cite{RK}, we define, for positive elements $a$ and $b$:
\begin{definition}  $a\geqCu b$ if there exists a sequence of elements
$(r_n)$  such that $r_n a r_n^{*}$ converges to $b$ in the norm topology. \label{def:geqCu}
\end{definition}\begin{definition}  If a positive element $a$ is nonzero and $a \geqCu a \oplus  a,$ then a is
said to be properly infinite.\end{definition}

We recall that a sufficient condition for the purely large property in the simple case is:
 \begin{proposition}\label{lem:characterize.purely.large}
Let $A$ and $B$ be separable \Cstar-algebras, with
$B$ stable.
Consider  a unital essential  extension of $B$ by~$A$. If the image of its Busby map $\tau$
 in the corona $\Mbb$ has the property that its positive nonzero elements are  properly infinite and full, then the extension is purely large.
\end{proposition}
\begin{proof}
Recall that the extension algebra $D$ is related to the given Busby map $\tau$ by a pullback construction \cite[\S7.2]{Loring}, where $\pi\colon\Mb\to\Mbb$ denotes the canonical quotient map:
\begin{center} \begin{tikzpicture}
    \node (D) at (0,0) {$D$};
    \node (M) at (2,0) {$\Mb$};
    \node (Q)  at (2,-1.6) {${\Mbb.}$};
    \node (A) at  (0,-1.6) {$A$};
    \draw[->] (D)--(A);
    \draw[right hook->] (A)--(Q) node[midway,above] {$\tau$};
    \draw[->] (D)--(M);
    \draw[->] (M)--(Q) node[midway,right] {$\pi$};
\end{tikzpicture} \end{center}
Regarding the extension algebra $D$ as a subalgebra of $\Mb,$ we see from the diagram that a positive element of $D$ that is not in the canonical ideal $B$ maps to a positive nonzero element of $\Mbb.$ Note that the essentiality of the extension means that the map $\tau$ is injective.
Therefore, if $c\in D$ is positive and not in the canonical ideal, then the nonzero element $\tau(c)$ is by hypothesis full and purely infinite. To show that the desired purely large property holds, we must show that the hereditary subalgebra
$$H:=\overline{cBc}$$
contains a stable full subalgebra. But we already noted that $\tau(c)=\pi(c)$ is properly infinite, and also full, in $\Mbb.$ By
Proposition 3.5 in \cite{RK} it follows that $\pi(c)\geqCu
1_{\Mbb},$ where $\geqCu$ denotes the Cuntz subequivalence relation of positive elements in $\Mbb,$ see \cref{def:geqCu}.
Choosing a
sufficiently large $n$ in  \cref{def:geqCu}, the operator $r_n\pi(c)r_n^{*}$ is then invertible.
Lifting to the multipliers, there is an $\tilde {r}\in\Mb$ such that
$\tilde {r}c\tilde {r}^{*}=1_{\Mb}+b$, where the element $b$ belongs
to the canonical ideal $B.$ Since $B$ is stable, there is a sequence of
isometries, $v_i,$  such that $v_i^{*}bv_i$ goes to zero in
norm\cite{HR}. We conclude that for some index $i$, the
expression $v_i^{*}\tilde {r}c\tilde {r}^{*}v_i=1+v_i^{*}bv_i$ is close enough to 1 to be invertible,
and thus there is an $r'\in \Mb$ such that $r'cr'{}^{*}=1_{\Mb)}$.
This implies that  $V:=c^{
1/2}r^{\prime *}$ is an infinite
isometry, so that its range projection $VV^*$ therefore has the purely large property. Since $VV^{*}\leq c\|r'\|^2$ ,  the hereditary subalgebra generated by
$VV^{*}$  is contained in the hereditary subalgebra generated by $c,$ and this establishes the purely large property for $c.$
\end{proof}
Recall that a representation is called \textit{essential} if its range as a map into $\Bh$  does not have any nonzero compact elements.
Kasparov considered  extensions of $B \otimes K$ by a separable  \C* $A,$ with Busby map induced by
\[
A \ \hookrightarrow \ 1 \otimes B ( \hs ) \ \hookrightarrow \ {\mathcal M} (B \otimes  K),
\]
where $A \hookrightarrow B(\mathcal H)$ is some faithful essential representation of
$A$ on the separable infinite-dimensional Hilbert space~$\hs$, and showed that it is, when not unital, absorbing in Kasparov's sense \cite[pg. 560, Lemma 1]{KasparovKK}. A very similar proof shows that, in our terminology, when unital, the extension is absorbing in the nuclear sense (i.e. it absorbs trivial unital and weakly nuclear extensions). This extension is furthermore trivial in the nuclear sense (lemma~12 in~ \cite{EK}).
We will call the Busby map of this extension the Kasparov extension, in honor of Kasparov's work, and will denote this map $\Kasp\colon A\to\Mbkbk$. The Kasparov extension is not unique, unless an equivalence relation is applied, but such an equivalence relation is implicit in Kasparov's theory. Hence we may as well refer to this extension as the Kasparov extension.
The main  property of  the Kasparov extension $\Kasp$ that we will use is that the range of the map $\Kasp$ is contained in a copy of the
Calkin algebra that is unitally embedded in the corona, $\Mbkbk$.

Recall Cuntz and Krieger's remark \cite[Remark 2.15]{cuntzkrieger} that the Cuntz algebra  $\O{\infty}$ can be defined concretely by isometries $v_i$ with the properties that $v_i^* v_j=0$
when $i\neq j$, and  $\sum v_iv_i^*=1$, *-strongly.
Proposition 23 in \cite{BasisProblem} constructs a copy of $\O{\infty}$ in a multiplier algebra and shows it can be chosen to have a similar property in both the multiplier algebra and the corona algebra.
Kirchberg~\cite[Remark 5.1]{KirchbergUnpublished} defines a unital *-homomorphism \label{def:Kirchberg}
$\delta_\infty\colon \Mb\rightarrow\Mbk$ given by
\[
\delta_{\infty}\colon m\mapsto\sum_{n=1}^{\infty}v_n^{}mv_n^{*}
\]
where the $v_n$ are isometries coming from this copy of the Cuntz algebra $\O{\infty}$ in the multipliers of $\Mbk$.
Kirchberg terms this map an infinite repeat. This terminology can be justified by observing that the above map is a sum
of *-homomorphisms of the form $h_i\colon m\to v_i m v_i^*$. These injective *-homomorphisms~$h_i$ have orthogonal ranges,
because the isometries are orthogonal, and thus the range of $\delta_{\infty}$ contains an infinite orthogonal sum of
copies of $\Mb$, contained within $\Mbk$. Each copy of $\Mb$ is embedded as a hereditary subalgebra in $\Mbk$ with unit $p_i:=h_i(1)$.
Each $h_i$ maps the simple essential ideal $B\subseteq\Mb$  into the simple essential ideal $\bk$.
Since each individual *-homomorphism $h_i\colon \Mb\to\Mbk$ is an isomorphism of $\Mb$ onto its range we can embed
$\Mbb$ inside the larger algebra $\Mbkbk$.
We will denote the component *-homomorphisms by
\[
h_i\colon\Mbb\to\Mbkbk,
\]
and their ranges will be called \textit{corona blocks\/}. The unit of a corona block is denoted~$p_i$.
Each corona block is an isomorphic copy of $\Mbb$.
Letting
\[
q_{\bk} \colon \Mbk\to\Mbkbk
\]
denote the canonical quotient map, the composition
\[\Sinf:=q_{\bk}\circ\delta_{\infty}\colon \Mb\to\Mbkbk\]
is the Busby map of an extension.
For technical reasons, usually it is desirable to restrict the domain to some separable, unital, and exact subalgebra of $\Mb$.
We call such a map a Kirchberg--Lin extension because of Lin's pioneering absorption result \cite[Theorem 1.12]{Lin}
for maps of this type. See also \cite{KirchbergUnpublished,EK}.
We summarize some known results in the following proposition. In order to simplify the language used, and in view of the fact that we will assume nuclearity, we say unitally absorbing instead of absorbing in the nuclear sense. It is understood that the extensions absorbed are unital.

\begin{proposition}\label{prop:abs}
Let $B$ be a separable simple \C*. Let $A$ be a separable and unital \Cs* of $\Mb$.
Then, if either $A$ or $B$ is nuclear, both the Kirchberg--Lin extension
$\Sinf\colon A\to \Mbkbk$ and the Kasparov extension $\Kasp\colon A\to \Mbkbk$ are unitally absorbing, unital, and trivial.
They are therefore unitarily equivalent, so that
\[
\Kasp(a)=U^*\Sinf(s)U \qquad (a\in A)
\]
for some unitary $U\in\Mb$.
\end{proposition}
\begin{proof}
That the Kirchberg--Lin extension is unitally absorbing when $A$ or $B$ is nuclear (which implies weak nuclearity) is shown in theorem 17.iii of \cite{EK},
see also \cite[Theorem 1.12]{Lin}. That the Kasparov extension is unitally absorbing is shown in   \cite{EK},
see also \cite{KasparovSV}. These unital extensions are trivial, as already discussed, and
unitally absorbing trivial unital extensions   are necessarily unitarily equivalent.
\end{proof}

We state a less technical corollary.
For a similar early result, with $D$ commutative, see~\cite[pg. 3030]{commutative}.

\begin{corollary}\label{cor:subalgebras}
Let $B$ be a separable, nuclear, and simple \C*. Suppose that $D$ is a norm-closed separable unital subalgebra, self-adjoint or not, of the range of the map
$\Sinf\colon \Mb\rightarrow\Mbkbk$.  Then there exists a unitary $U\in\Mb$ such that $U^* D U$ is contained
in a copy of the Calkin algebra within $\Mbkbk$.
\end{corollary}
\begin{proof}
If $D$ is not self-adjoint, let $C^*(D)$ denote the (separable) \C* it generates.
Since $\Sinf$ is injective, let $A:=(\Sinf)^{-1}(C^*(D)).$ This is clearly  a unital separable \Cs* of $\Mb.$
The Kirchberg--Lin extension
$\Sinf \colon A\rightarrow\Mbkbk$  is unitarily equivalent to the Kasparov extension $\Kasp\colon A\to \Mbkbk$.
We recall that the range of the Kasparov extension is contained in a unitally embedded copy of the Calkin algebra within
$\Mbkbk$, and thus the unitary equivalence of the extensions implements the desired equivalence of subalgebras.
\end{proof}
\begin{remark}\label{rem:inf.repeats}
The above result already implies, for example, that nonzero elements of the form $\Sinf(m)$
cannot be contained in any proper ideal of the corona --- this is because the Calkin algebra is simple
and unitally embedded, so that its nonzero elements are not contained in any proper ideal of the corona.
Thus we have a short proof of a slight generalization of a
 known result \cite{brown,Lin} that constant infinite repeats are strongly full.
\end{remark}


\section{On the structure of double relative commutants}
\noindent
Suppose that $A$ is some given unital separable subalgebra of a multiplier algebra $\Mb$, with $B$ simple and separable.
The Busby map of the associated Kirchberg--Lin extension is the map $\Sinf\colon A\to\Mbkbk$,
given by the composition of Kirchberg's homomorphism $\delta_\infty$ with the canonical quotient by $B\otimes \mathcal K$.
Letting $q_{\bk} \colon \Mbk\to\Mbkbk$ and $q_B \colon \Mb\to\Mbb$ denote the canonical quotient maps, we note that, for all~$a$,
\begin{equation}
h_i\circ q_B (a)=(q_\bk \circ h_i)(a)=q_\bk (p_i\constant{a} p_i)=p_i q_{\bk}(\constant{a})p_i.
\end{equation}
This identity shows that the subalgebra $\Sinf(A)\subseteq\Mbkbk$ contains copies of $q_B(A)$
unitally contained in hereditary subalgebras $\her(p_i)$, where each copy comes from the map $h_i$ defined earlier. We will denote the image of $A$ in the corona block $\her(p_i)$ by $D_i$,
and moreover we note that each subalgebra $D_i$ is contained in an isomorphic copy of $\Mbb$, namely,
the range of the map $h_i\colon\Mbb\to\Mbkbk$.
We now consider the properties of the subalgebra $D_i.$

\begin{lemma}\label{lem:relative}
Let $B$ be a simple separable \C*. Let $p$ be a projection in the corona
$W:=\Mbb$ of this algebra. Let $D$ be a unital subalgebra of $p W p$,  {self-adjoint or not}. Then the  double relative commutant of $D$ is the same,
up to a unitisation, whether it is relative to $pWp$ or to~$W$.
\end{lemma}
\begin{proof}
Let $L:=pWp$. First we will show that an element $x$ of $D_{W}''$, the double relative commutant of $D$ in $W$,
can be decomposed relative to $p$ as
$\begin{pmatrix} \lambda\,(1-p) & 0 \\
                 0 & pxp
                 \end{pmatrix} $ where $\lambda$ is a scalar.
Then we will show that $pD_{W}''p = D_{L}''$, and this will prove the result.

Since $D$ is contained in $L$, the double relative commutant $D_{W} ''$ of $D$ relative to $W$ is contained in $L_{W}''$.
 \cref{th:main} shows that the double relative commutant of the \Cs* $L$ in $W$ is $L$ unitised by the unit of the corona,
and thus $D_{W}''$ is contained in the unitisation of  $L:=pWp$. In other words, the elements of $D_{W}''$ can be decomposed as shown above.

Next, we show that $pD_{W}''p = D_{L}''$, by proving two inclusions.
We begin by showing that $D''_{L}$ is contained in $D''_{W}$.
Consider an element $s$ of the relative commutant $D'_W$.
Since the given projection $p$ is the unit of $L$, and $D$ is a unital subalgebra of $L$,
it follows that $p$ is in $D$. Therefore, $s$ commutes with $p$.
Thus $s$ is diagonal with respect to $p$, meaning that $s=psp + (1-p)s(1-p)$.
The first term $psp$ is in $D'_L$, and the second term $(1-p)s(1-p)$ is in the annihilator of $L$, namely $\her(1-p)$.
Consequently $D'_W\subseteq D'_L + \her(1-p)$.
On the other hand, it is clear that $D'_L$ and $\her(1-p)$ are both in $D'_W$, so  we have the reverse inclusion as well.
Therefore,
\[
D'_W = D'_L + \her(1-p).
\]
Notice that an element of  $D''_{L}$ will commute with elements of $D'_L$ and will annihilate elements of $\her(1-p)$.
Thus, such an element will commute with elements of the above right hand side, and therefore commutes with $D'_W$.
This proves that $D''_{L}$ is contained in $D''_{W}$.

Finally, we show that $pD_{W}''p \subseteq D_{L}''$. Since the left hand side is clearly contained in $L$, we need only check
that the elements on the left hand side commute with~$D'_L$. However, $D'_L$ is a subset of $D'_W$ and  $D_{W}''$
commutes elementwise with $D'_W$ --- implying that $D''_W$ commutes elementwise with $D'_L$.
The element $p$ acts on $D'_L$ as the unit, thus also commutes with $D'_L$. Therefore we have the required inclusion.
\end{proof}

\newcommand{\F}{\mathcal F}
\def \V #1#2{v_{#1,#2}}
\def \Vs #1#2{v^{*}_{#1,#2}}
\def \P #1{p_{#1}}
\def \U #1#2{u_{#1,#2}}
\def \Us #1#2{u^{*}_{#1,#2}}
We now define some convenient elements that will be used within \cref{lem:diagonals,theorem:RangeOfExtensionDCT}.
Let $\{\P n \}_{n=1}^{\infty}$ be a sequence of pairwise orthogonal projections in the multipliers $\Mb$ of a stable \Cstar-algebra $B$ such that
$$\P n \sim 1 \mbox{ for all $n$}$$
and
$$\sum_{n=1}^\infty \P{n}=1$$ with strict convergence (see \cite{HR}, \cite[Prop. 23]{BasisProblem}). The $\P n$ are elements of $1\tensor\Bh\subset\Mb.$ For all $n\geq1,$ let $\V n 1 \in1\tensor\Bh\subset\Mb$
be a partial isometry such that
$$\Vs n 1 \V n 1 =\P{1} \mbox{ and } \V n 1 \Vs n 1 = \P n $$
For all $m,n,$ we define $$ \V n m :=\V n 1 \Vs m 1. $$
Finally, let $\F$ be the collection of unitaries in $\Mb$ containing all unitaries interchanging the $p_n$ pairwise. We suppose that $\F$ contains all unitaries of the form $\U mn = \V mn + \V nm + (1-\P m - \P n ),$ and furthermore  is closed under adjoints, multiplication, and the strict topology. The strict topology restricted to $1\tensor\Bh$  coincides with a strong topology.

 \begin{lemma}\label{lem:diagonals}
Let $A$ be a unital sub-\Cstar-algebra of $\Mb,$ and let $B$ be stable. Then in the relative commutant $\Sinf(A)'$ there is a unitally embedded copy of $\O{n},$ for $n=2,3,\cdots,\infty.$ The relative commutant $\Sinf(A)'$ contains the unitaries $\F.$
\end{lemma}
\begin{proof}
The algebra $\delta_\infty(A)$ consists of elements that are diagonal with respect to the projections $p_n$.
Since the generators $\U mn$ of $\F$ satisfy $\delta_\infty(A)=\U mn ^* \delta_\infty(A)\U mn$ it follows that $\F$ is in the relative commutant $\Sinf(A)',$ as claimed.
These unitaries $\U mn$ generate a unitally embedded copy of $B(\hs)$ in the multiplier algebra, see for example Lemma 5.2.ii in \cite{KirchbergUnpublished}.
{Alternatively, one can proceed as in the proof of lemma 1 in \cite{paschke}}. 
We can find in this copy of $B(\hs)$  a unitally embedded copy of the Cuntz algebra $\O{2}.$
We thus have  a unitally embedded copy of $\O{2}$ in the relative commutant of $\Sinf(A)'$ in the corona. The case of $\O{n}$ is similar. 	
\end{proof}

\begin{lemma}\label{lem:kk.diagonals}
Let $A$ be a unital separable sub-\Cstar-algebra of $\Mb$, and let $B$ be simple, stable, nuclear, and separable.
If $S$ is a unital, exact,  and separable subalgebra of $\Sinf(A)''$, then $S$ is contained in the range of a trivial absorbing extension.
The positive elements of $S$ are full and properly infinite.
\end{lemma}
\begin{proof}
The given subalgebra $S$ commutes with $\Sinf(A)'$, and by
\cref{lem:diagonals}  it follows that $S$  commutes with a unital copy of $\O{n}, n=2,\cdots,\infty$ that comes from the multiplier algebra.
Denoting the generators of the copy of $\O{2}$ by $w_i$,
the unital inclusion map $\tau\colon S\rightarrow \Mbb$ then has the
property that $\tau(s)= w_1 \tau(s) w_1^*+ w_2 \tau(s) w_2^*$ for all $s\in S.$ This implies that positive elements of the form $\tau(s)$ are purely infinite and fullness is similar, as in  \cref{rem:inf.repeats}.
Since the copy of $\O{2}$ comes from the multiplier algebra, we have here exactly the definition of BDF addition, \textit{i.e.,} addition of extensions. But this implies that as an extension, $\tau=\tau+\tau$ in the enveloping abelian group of extensions. In a group, the only element satisfying $\tau+\tau=\tau$ is the trivial element. Thus the extension $\tau$
is  trivial in the enveloping abelian group of extensions.
 Since $\tau$ is a unitally absorbing extension, by \cref{lem:characterize.purely.large,lem:inf.repeats}, triviality in the group implies being unitarily equivalent by a multiplier unitary to, for example, Kasparov's unitally absorbing trivial (\textit{i.e.,} split) extension. But then the extension $\tau$  splits as well, so is trivial as claimed.
\end{proof}

Since any positive element $x$ generates an abelian unital nuclear separable subalgebra $\Cstar(x,1)=:S,$  we can deduce from the above that every positive element of $\Sinf(A)''$ is purely infinite and full:
\begin{corollary}\label{cor:range.is.pi}
Let $B$ be simple, stable, nuclear, and separable. Let $A$ be some unital and separable subalgebra of $\Mb.$ The positive elements of $\Sinf(A)''$ are full and properly infinite.
\end{corollary}

The above result is the key step needed in the next section.
 \newcommand{\Sphi}{\bar\delta}
We also mention that, as can be shown by a direct method, the elements of $\Sinf(A)''$ are contained in the range of $\Sinf,$  and this again implies that the nonzero positive elements of $\Sinf(A)''$ are full and properly infinite:
\begin{theorem} Let $A$ be a unital separable subalgebra of $\Mb$, and let $B$ be simple, stable, nuclear, and separable. The elements of $\Sinf(A)''$ are contained in the range of $\Sinf.$
\label{theorem:dd.is.diagonal} \label{theorem:RangeOfExtensionDCT}
\end{theorem}
\begin{proof}
\def \A{a_0}
Choose some element $a$ of $\Sphi(A)''$ in $\Mbb.$ Lift this element to $\A$ in $\Mb,$ and let $\A':= \A-\sum \P n \A \P n,$ where the sum converges strictly.
We  show that $\A'$ is actually an element of $B.$ We are free to make the usual small adjustments to a lifting, in particular, we are free to conjugate by elements of $\mathcal F$ because they commute
 with the given element in the corona. The strategy will be to make such adjustments and then to show that $\P {m,n} \A' \P {m,n},$ which is compact in the Hilbert module sense, where $\P {m,n}$ denotes $\P m + \P {m+1} +\cdots \P n,$ has the Cauchy property with respect to $m$ and $n,$ and this will show that the limit, namely $\A',$ is compact in the Hilbert module sense, or in other words, is in $B.$

Recall that $\V n 1$ denotes a multiplier partial isometry such that $\Vs n1 \V n1 = \P 1 ,$  $\V n1 \Vs n1 = \P n ,$ and $\V nm = \V n1 \Vs m1.$
Since $\mathcal F$ commutes with $\Sphi(A),$ we observe:
\begin{enumerate}\item[(i)] $\P m \A \P n\in B$ when $n$ and $m$ are distinct,
\item[(ii)] $u\A u^* - \A'\in B$ for all $u\in\mathcal F,$
\item[(iii)] $u\A'u^* - \A'\in B$ for all $u\in\mathcal F,$
\item[(iv)] $\P m \A \P m - \V mn \P n \A \P n \V nm\in B$ for all $m, n,$
\item[(v)] Given an $\varepsilon>0,$ if $m,n\geq N(\varepsilon)$, then $|| \P m \A \P m - \V mn \P n \A \P n\V nm ||<\varepsilon,$
\item[(vi)] $(\P 1 + \P 2 + \cdots + \P n )\A'\in B$ for all $n,$
\item[(vii)] $[\P n,\A]\in B$ for all $n,$ and
\item[(viii)] $[\P n,\A]$ goes to zero in norm as $n$ goes to infinity.
\end{enumerate}
Thus, for all $1\leq m\leq n,$ for all $\varepsilon>0,$ and for all $\ell\geq 1,$ there exist
integers $\ell_j$ with $\ell<\ell_1<\ell_2<\cdots<\ell_{n-m+1}$ such that
\[
\|  \P {m,n} \U m{\ell_{n-m+1}} \cdots \U {n-1}{\ell_2} \U n{\ell_1}
                                           \A' \Us n{\ell_1} \Us {n-1}{\ell_2} \cdots \Us m{\ell_{n-m+1}} \P {m,n}  \|<\varepsilon
\]  Thus $\A'$ is in $\bk,$ or \[\pi(\A)=\pi(\sum_{n=1}^\infty \P n \A \P n)\]
where $\pi$ denotes the canonical quotient map.
It remains to show that the right hand side is in the range of $\Sphi.$ For this, consider
the norm limit $d:=\lim_{n\rightarrow\infty} v_n^* p_n a_0 p_n v_n .$ Then, the element $d$ is in $\P 1 \Mb \P 1$ and  $\sum_{n=1}^\infty \V 1n d \Vs 1n $ is an operator in $\Mb$ which is a lift of $\A .$
\end{proof}

\section{Arveson's distance formula}
\noindent

The original form of Arveson's distance formula  applies to norm-closed subalgebras of the classic Calkin algebra, self-adjoint or not,  and can be phrased as follows:
\begin{lemma}[\!\!{\cite[p. 344]{arveson}}]
Consider the Calkin algebra of a separable infinite Hilbert space. If $D$ is a separable unital norm-closed subalgebra of the  Calkin algebra, and $x$ is an element of the Calkin algebra, then there exists a projection $p$ such that $p$ commutes with elements of $D$
and
\[
\dist(x,D)=\|(1-p)x p\|.
\] \label{lem:originalArveson}
\end{lemma}

The following lemma is similar to \cite{GN}, lemma~3.3. It applies to, for example, separable subalgebras of the range of the trivial extension $\Sinf.$
\begin{lemma}\label{lem:fixed.gn}
Let $B$ be a separable, stabilized, and nuclear \C*. Let $D$ be a separable norm-closed unital subalgebra, {self-adjoint or not},
 of the range of any trivial unital extension. Suppose that the   nonzero positive elements of the range are full and properly infinite.
and let $x$ be an element in that range. Then there exists a projection $p\in \Mbb$ such that $p$ commutes with elements of $D$
and
\[
\dist(x,D)=\|(1-p)x p\|.
\]
\end{lemma}
\begin{proof}
Suppose that the given trivial extension is actually $\Sinf,$ and consider two cases.
Suppose first that $B$ is isomorphic to the compact operators.
Then, $\Mbb$ is isomorphic to the Calkin algebra, and we apply \cref{lem:originalArveson} to the subalgebra $D$ and the element~$x$.

Now, suppose that $B$ is not isomorphic to the compact operators.
Then  there is a unitally embedded copy of the Calkin algebra within $\Mbb$, and we use \cref{cor:subalgebras}
to find a unitary such that $UDU^*$ and  $UxU^*$ are in this copy of the Calkin algebra.
We then apply the previous case to find a projection such that $\dist(UDU^*,UxU^*)=\|(1-p)UxU^*p\|$.
But both the distance and the norm are unitarily invariant, so it follows that $\dist(D,x)=\|U^*(1-p)U x U^*pU \|$.
The projection $U^*pU$ therefore has the required properties. This proves the lemma in the case where the trivial extension is $\Sinf.$
In the general case, the given trivial extension $\tau$ is  unitarily equivalent to $\Sinf$, and we have seen that unitary equivalence is sufficient.

\end{proof}

\section{Main results}

\begin{theorem}\label{theorem:main}
Let $B$ be a separable simple stable nuclear \C*.
Suppose that $A$ is a separable norm-closed unital subalgebra of $\Mb$. Then $\Sinf(A)$ is equal to its double relative commutant.
\end{theorem}
\begin{proof}

We are to show that  $\bar x \in \Sinf(A)''$ is contained in $\Sinf(A).$
Then $\bar x$ and $\Sinf (A)$ are contained in a unital separable \Cs* $E$ of $\Sinf(A)''$. This algebra $E$ might not be nuclear, but
by \cref{cor:range.is.pi}, the  nonzero positive elements of $E$ are properly infinite and full, and  by \cref{theorem:RangeOfExtensionDCT}, the \Cs* $E$ is contained in the range of a trivial extension.

Then this is sufficient to apply the distance formula  of \cref{lem:fixed.gn}.
Thus, we have a  projection $p\in\Mbb$ such that $p$ commutes with elements of $\Sinf(A)$ and
$\dist(\bar x,\Sinf(A))=\|(1-p)\bar x p\|$. Since $\bar x$ must commute with $p$,
the element on the right is zero. Thus
$\dist(\bar x,\Sinf(A))$ is zero, which implies that $\bar x $ was actually in~$\Sinf(A)$, as was to be shown.
\end{proof}


\begin{lemma} The projections $p_i=v_i v_i^*$ from page \pageref{def:Kirchberg} sum strictly to 1 in the multipliers. The supremum of the finite sums of the $p_i$ in the corona is 1.\label{fatalni.lemma}
\end{lemma}\begin{proof} 

For the strict convergence, see \cite[p.155]{HR}. This implies that the supremum of the finite sums of the $p_i$ in the multiplier algebra is 1. The quotient map  into the corona is surjective, and surjective maps preserve suprema, so the supremum in the corona is $1.$
For a related discussion and supplementary information, see \cite{BasisProblem}, especially Proposition 23.
\end{proof}

\begin{lemma}\label{lem:technical.ii}
Let $p=h_i(1)$. Then
\begin{equation}\label{eq:commutant}
p \Sinf(A)' p = (p\Sinf(A)p)',
\end{equation}
where the commutant on the right is relative to $\her(p)$, and the commutant on the left is relative to the whole corona.
\end{lemma}
\begin{proof}
Since $p$ is in $\Sinf(A)'$, $p\Sinf(A)'p$ is a subalgebra of $\Sinf(A)'$ and thus commutes with
$\Sinf(A)$. Now $p\Sinf(A)'p$ trivially commutes with $p$, and therefore commutes also with $p\Sinf(A)p$. Thus,
\begin{equation*}
p \Sinf(A)' p \subseteq \left( p\Sinf(A) p \right)'.
\end{equation*}
For the reverse inclusion, let $y$ be an element of $(p \Sinf(A) p)'$
relative to the  corona block with unit~$p$. Since $y$ is contained in $\her(p)$ it commutes with $(1-p)\delta(A)$. Therefore $y$ is in $\Sinf(A)'$ and $y$ is in $\her(p)$.
This means that $y$ is in $p \Sinf(A)' p$ as claimed, and this establishes identity~\eqref{eq:commutant} above.
\end{proof}

\begin{proposition}\label{prop:technical.i} Let $B$ be stable.
The double relative commutant of\/ $D_i$ relative to its corona block, $\her(p)$,
is contained in $p\Sinf(A)''p$, where $\Sinf(A)''$ is relative to the whole corona.
\end{proposition}
\begin{proof}
Suppose that $w$ is an element of the double relative commutant of\/ $D_i$ relative to its corona block, $\her(p)$. There exists some element $m\in\Mb$ such that $w=p\Sinf(m)p.$ Then $p_j\Sinf(m)p_j$ is, for each $j,$ unitarily equivalent to $w$ by one of the unitaries interchanging corona blocks provided by \cref{lem:diagonals}. By unitary equivalence,  $p_j\Sinf(m)p_j$ not only is in the corona block $\her(p_j)$ but is in $D_j''$. We are to show that $\Sinf(m)$ commutes with $\Sinf(A)'.$

 \cref{lem:technical.ii,lem:diagonals} show that  $p_i\Sinf(A)'p_i = D'_i$ and that  $p_i\Sinf(A)'p_j= p_i D_i' U_{ij}p_j,$ where $U_{ij}$ is the unitary that intertwines $p_j$ and $p_i.$
Then, $p_i\Sinf(A)'p_j$ and $\Sinf(m)$ commute.  This means that for each $a$ in $\Sinf(A)',$ we have an element $z=[a,\Sinf(m)]$ having the property that $p_i z p_j$ is zero for all $i$ and $j$. 
By \cref{fatalni.lemma} it follows that if $p_i z p_j$ is zero for all $i$ and $j$  then
 $z$ is zero. But then $\Sinf(m)$ is in $\Sinf(A)'',$ as was to be shown.
\end{proof}
Applying $p_i$ from left and right to the result of  \cref{th:main}, we have a corollary.

\begin{corollary}\label{cor:nonstable}
Let $B$ be a separable simple nuclear \C*.
Suppose that $D$ is a separable norm-closed unital subalgebra of $\Mbb$. It is then equal to its double relative commutant.
\end{corollary}
\begin{proof} If $B$ is stable, the result follows from \cref{theorem:main} and \cref{prop:technical.i}. If $B$ is not stable, then the given unital subalgebra $D$ of the corona $\Mbb$ is a subalgebra of the $p_1$ corner of the corona of the stabilization,
$\Mbkbk.$ The previous case implies that the double relative commutant of $D$ relative to the corona of the stabilization is equal to the unitisation of~$D$.
 Now \cref{lem:relative}, with $p$ taken to be the unit of the corner, implies that the double relative commutant of $D$ remains the same, up to a unitization, when computed in $p(\Mbkbk)p = \Mbb$.
\end{proof}

\section*{Ackowledgements}
Both authors were supported by the Fredrik and Catherine Eaton Fellowship Fund
while the first-named author was also supported by NSERC (Canada).

\end{document}